\documentclass[a4paper, 12pt, reqno]{amsart}
\usepackage{amsmath, amssymb, latexsym, graphicx, hyperref}

\newcommand{\R}{\mathbb{R}}

\newcommand{\w}{\omega}
\newcommand{\worder}{$\omega$-or\-der}
\newcommand{\wordered}{$\omega$-or\-de\-red}

\newcommand{\wpartition}{$\omega$-par\-ti\-tion}
\newcommand{\wpartitions}{$\omega$-par\-ti\-tions}
\newcommand{\wpartitioned}{$\omega$-par\-ti\-tio\-ned}

\newcommand{\wasymmetrical}{$\omega$-a\-sym\-me\-tri\-cal}

\begin{document}

\title{Transfinite partitions of Jordan Curves}

\maketitle

\begin{center}
\vspace{-16pt} \small{ Antonio Leon Sanchez
(\href{mailto:aleons@educa.jcyl.es}{aleons@educa.jcyl.es})\\I.E.S.
Francisco Salinas. Salamanca, Spain.}
\end{center}

\pagestyle{myheadings}

\markboth{\small{Transfinite partitions of Jordan
Curves}}{\small{Transfinite partitions of Jordan Curves}}

\begin{abstract}
The $\w$-asymmetry induced by transfinite partitions makes it
impossible for Jordan curves to have an infinite length.
\end{abstract}

\section{$\w$-asymmetry}

\noindent As we known from the XVIII century, \wpartitions\ (as we
call them nowadays) of finite line segments are only possible if
the successive adjacent parts of the \wpartition\ are of a
decreasing length. This inevitable restriction induces a huge
asymmetry in the very partition. In fact, whatever be the length
of the \wpartitioned\ line segment and whatever be the
\wpartition, all its parts, except a finite number of them, will
necessarily lie within an arbitrarily small final segment. For the
sake of illustration, consider an \wpartition\ of a $10^{30}$
light years length segment -the assumed diameter of the universe.
Whatever be the \wpartition\ of this enormous line segment all its
infinitely many parts, except a finite number of them, will
inevitably lie within a final segment inconceivable less than, for
instance, Planck length ($\sim 10^{-33}$ cm). There is no way of
performing a more equitable partition if the partition has to be
\wordered. Thus, \wpartitions\ are \wasymmetrical. For the same
reason it is impossible to consider two proper points in the real
line $\R$ separated by an infinite euclidean distance, in spite of
the assumed infiniteness of the real line. The above simply
unaesthetic consequences of \wpartitions\ become a little more
controversial if the partitioned object is a closed line as a
Jordan curve. The objective of the following short discussion is
just to examine one of those consequences.

\begin{figure}[htb!]
\centering
\includegraphics[scale=0.9]{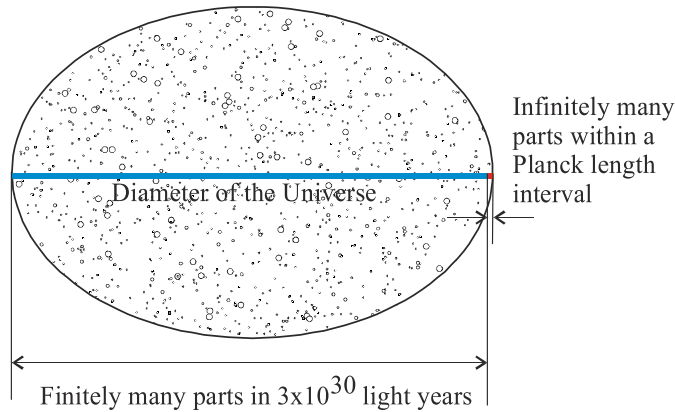}
\caption{A cosmic $\w$-asymmetry. The \worder\ makes it impossible
a more equitable distribution of the available space.}
\end{figure}

\section{Transfinite partitions of Jordan Curves}

\noindent Let $f(x)$ be a real valued function whose graph is a
Jordan Curve $\mathbf{J}$ in the euclidean plane $\mathbb{R}^2$.
If $a$ and $b$ are any two $\mathbf{J}$'s points, we will write
$L(a,b)$ to denote the length of the $\mathbf{J}$'s arc
$\widetilde{ab}$ whose endpoints are $a$ and $b$. That is to say:

\begin{equation}
L(a,b) = \int_a^b \sqrt{1+(f(x)')^2}dx
\end{equation}

\noindent \\Assume that $\mathbf{J}$ has an infinite euclidean
length. In these conditions let $r$ be any proper real number
greater than 0 and assume $\mathbf{J}$ is partitioned clockwise
from a point $x_1$ into a certain number of adjacent parts
$\widetilde{x_1x}_2$, $\widetilde{x_2x}_3$, $\widetilde{x_3x}_4$
$\dots$ so that each part $\widetilde{x_ix}_{i+1}$ has a finite
length which is equal or greater than $r$:

\begin{equation}\label{eqn:FJC >= r}
L(x_i, x_{i+1}) \geq r, \ \forall i \in I
\end{equation}

\noindent \\Evidently these partition $\langle x_i \rangle_{i \in
I}$ must be infinite otherwise, and being finite the length of the
parts, $\mathbf{J}$ would have a finite length. In addition, and
according to Cantor \cite{Cantor1966}, $\langle x_i \rangle_{i \in
I}$ cannot be uncountably infinite. In fact, consider the sequence
of real numbers $\langle r_i \rangle_{i \in I}$ defined as:

\begin{align}
&r_1 = x_1\\
&r_{i+1} = r_i + L(x_i, x_{i+1}), \, \forall i \in I
\end{align}

\noindent \\The one to one correspondence $f$ between $\langle x_i
\rangle_{i \in I}$ and $\langle r_i \rangle_{i \in I}$ defined by
$f(x_i) = r_i$ proves that both sequence have the same
cardinality. So if the first one were uncountably infinite so
would be the second. But $\langle r_i \rangle_{i \in I}$ cannot be
uncountably infinite because if that were the case we could pick a
different rational number $q_i$ in each real interval $[r_i,
r_{i+1})$ and then we would have an uncountable set of different
rational numbers, which is obviously impossible \cite{Cantor1966}.
The partition $\langle x_i \rangle_{i \in I}$ must therefore be
countably infinite. On the other hand, the adjacency of
$\widetilde{x_1x}_2$, $\widetilde{x_2x}_3$, $\widetilde{x_3x}_4$
$\dots$ implies successiveness and then ordinality. Accordingly
$\langle x_i \rangle_{i \in I}$ must have a certain ordinal
number, and taking into account its countable nature, this ordinal
can only be transfinite, i.e. an ordinal of the second class
according to Cantor's classic terminology \cite{Cantor1955}.

\begin{figure}[htb!]
\centering
\includegraphics[scale=1]{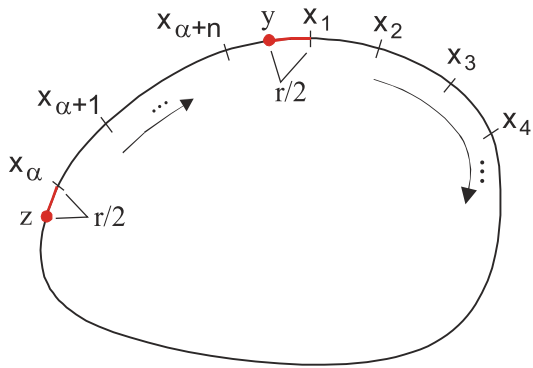}
\caption{Transfinite partition of a Jordan Curve in the euclidean
$\R^2$}
\end{figure}

Now consider a point $y$ anticlockwise from $x_1$ and such that
$L(y,x_1) = r/2$. According to (\ref{eqn:FJC >= r}) $y$ can only
belong to the last element of $\langle x_i \rangle_{i \in I}$.
Consequently this sequence must have a last element and then must
be $(\alpha+n)$-ordered, being $\alpha$ a transfinite ordinal of
the second class second kind\footnote{Ordinals of the second class
first kind have immediate predecessor and immediate successor,
while those of the second kind have immediate successor but not
immediate predecessor.}, and $n$ a finite ordinal; so that
$(\alpha + n)$ is an ordinal of the second class first kind. If
must therefore exist an $\alpha$-th second class second kind
ordinal in $\langle x_i \rangle_{i \in I}$ so that
$\widetilde{x_\alpha x}_{\alpha+1}$ cannot have an immediate
predecessor ($x_\alpha$ is the limit of a sequence of first kind
ordinals). Now then, according again to (\ref{eqn:FJC
>= r}), the point $z$ anticlockwise from $x_\alpha$ and such that
$L(z, x_\alpha) = r/2$ can only belong to a part immediately
preceding $\widetilde{x_\alpha x}_{\alpha + 1}$, which is
impossible if $x_\alpha$ is of the second kind, as it has to be.
This proves that $\langle x_i \rangle_{i \in I}$ must be, but
cannot be, $(\alpha+n)$-ordered, a contradiction we have derived
from our initial assumption on $\mathbf{J}$ infiniteness. We must
therefore conclude that Jordan Curves of infinite length are
inconsistent objects, and so must be any metrical space compatible
with them.

In addition, if $\mathbf{J}$ is finite the above condition
(\ref{eqn:FJC >= r}) can no longer holds, otherwise $\mathbf{J}$
would have an infinite length. Therefore, transfinite partitions
of Jordan curves can only be of a decreasing part-length, i.e.
$\w$-asymmetrical.


\begin{thebibliography}{1}

\bibitem{Cantor1955}
Georg Cantor, \emph{Contributions to the founding of the theory of
transfinite
  numbers}, Dover, New York, 1955.

\bibitem{Cantor1966}
\bysame, \emph{{\"{U}}ber unendliche lineare
{P}unktmannigfaltigkeiten},
  Abhandlungen mathematischen und philosophischen {I}nahalts (E.~Zermelo, ed.),
  Olms, Hildesheim, 1966, pp.~149 --157.

\end{thebibliography}

\providecommand{\bysame}{\leavevmode\hbox
to3em{\hrulefill}\thinspace}
\providecommand{\MR}{\relax\ifhmode\unskip\space\fi MR }
\providecommand{\MRhref}[2]{%
  \href{http://www.ams.org/mathscinet-getitem?mr=#1}{#2}
} \providecommand{\href}[2]{#2}

\end{document}